\documentclass[12pt]{article}
\usepackage{mathrsfs}
\usepackage{amsmath}
\usepackage{amssymb}
\usepackage{amsthm}
\usepackage{amsfonts}
\usepackage{amscd}
\usepackage[mathscr]{eucal}
\newcommand{\Z} {{\mathbb  Z}}
\newcommand{\Q}{{\mathbb  Q}}

\newcommand{\C}{{\mathbb  C}}

\begin{document}
\parindent  25pt
\baselineskip  10mm
\textwidth  15cm    \textheight  23cm \evensidemargin -0.06cm
\oddsidemargin -0.01cm

\title{ {On cohomology groups $ H^{1} $ of $ G-$modules of
finite type over cyclic groups }}
\author{\mbox{}
{ Derong Qiu }
\thanks{ \quad E-mail:
derong@mail.cnu.edu.cn, \ derongqiu@gmail.com } \\
(School of Mathematical Sciences,
 Capital Normal University, \\
Beijing 100048, P.R.China )  }

\date{}
\maketitle
\parindent  24pt
\baselineskip  10mm
\parskip  0pt

\par   \vskip 0.4cm

{\bf Abstract} \ Let $ G $ be a cyclic group, in this paper,
we study the Herbrand quotient and $ 1-$th cohomology group
on finitely generated $ G-$modules in some cases. When $ G $ is
of order $ 2, $ the order of the cohomology group is explicitly related to
some invariants, and this relation is used to study unit groups
over quadratic extensions of number fields. We also give some applications
on Pell equations and class number of number fields.
\par  \vskip  0.2 cm

\par  \vskip  0.1 cm


\par     \vskip  0.4 cm

\hspace{-0.8cm}{\bf 1. Introduction}

\par \vskip 0.2 cm

The $ 1-$th cohomology group $ \text{H}^{1} (G, A) $ and the
Herbrand quotient $ h(G, A) $ for a $
G-$module $ A $ are very useful for studying arithmetic when $ A $ is
an arithmetic object, e.g., $ A $ is related to elliptic curves or
number fields. However, their calculation are usually not easy, even
when $ G $ is a group of order $ 2. $

Let $ G $ be a cyclic group of order $ n $ with a generator $ \sigma , $
and let $ (A, + ) $ be a finitely generated  abelian group, we
denote $ r(A) = \text{rank}(A). $ We assume that $ A $ is a $
G-$module. Throughout this paper, for a set $ S, $ we denote its
cardinal by $ \sharp S. $ For an arbitrary abelian group $ B $ and
a positive integer $ m, $ we denote $ m B = \{m b : \ b \in B \} $ and
$ B [m] = \{b \in B : \ m b = 0 \}. $ For a $ G-$module $ A, $
one has the following Tate cohomology groups:
$ \widehat{\text{H}}^{i}(G, A) = \text{H}^{i}(G, A) \quad \text{if} \
i \geq 1; \quad \ \widehat{\text{H}}^{0}(G, A) = A^{G} / N_{G} A, $ \
where $ N_{G} = \sum _{\tau \in G} \tau \in \Z[G], $ and $ \Z[G] $
is the group algebra of $ G $ over $ \Z. $ \\
Since $ G $ is a finite cyclic group, $ \widehat{\text{H}}^{-1}(G, A) =
\text{H}^{1}(G, A). $ The Herbrand quotient of the $ G-$module $ A $ is
$ h(G, A) = \frac{\sharp \widehat{\text{H}}^{0}(G, A)}
{\sharp \text{H}^{1}(G, A)} $ (see [AW] and [Se]).

For the Herbrand quotient and the $ 1-$th cohomology group of a $ G-$module
$ A, $ we have
\par \vskip 0.2 cm

{\bf Proposition 1.} \ Let $ G $ be a finite cyclic
group, and $ A $ be a $ G-$module. If there exists a $ G-$submodule $ B $
of finite index in $ A $ which is a finitely generated free abelian group
containing a $ G-$invariant $ \Z-$basis $ X = \{x_{1}, \cdots, x_{r} \} $
(i.e., $ \tau X \subset X $ for all $ \tau \in G, $ or in other words,
$ X $ is a $ G-$set). Then
$$ h(G, A) = \prod _{x \in X / G} \sharp G_{x}, \quad \text{and} \quad
\sharp \text{H}^{1}(G, A) = \frac{(A^{G} : N_{G}A)}
{\prod _{x \in X / G} \sharp G_{x}}, $$ where $ G_{x} = \{\tau \in G :
\ \tau x = x \} $ is the stabilizer of $ x. $

The proof will be given in Section 2 (see Proposition 2.1 in the following).
\par \vskip 0.2 cm

When $ G = < \sigma > $ is of order $ 2, $ for a $ G-$module $ A $
which is a finitely generated abelian group, we simply denote
$ NA = (1+\sigma )A = \{a + \sigma a : a
\in A \}, \ N^{-}A = (1-\sigma )A = \{a - \sigma a : a \in A \}, \
A^{+} = A^{G}, \ r_{+}(A) = \text{rank}(A^{+}), $ and $ r_{-}(A) =
\text{rank}(A^{-}), $ where $ A^{-} = \{ a \in A : \ \sigma a = -a
\}. $ Obviously, $ r(A) = r_{+}(A) + r_{-}(A), \ A^{+} [2] = A^{+}
\cap A^{-} = A^{-}[2], \ (A : A^{+} + A^{-}) = (N^{-}A : 2A^{-}) =
(NA : 2 A^{+}). $
\par \vskip 0.2 cm

{\bf Proposition 2.} \ Let $ G $ be a group of order $ 2, $ and $ A $ be
a $ G-$module. If $ A $ is a finitely generated abelian group. Then
\begin{align*} \sharp \text{H}^{1} (G, A) &= \frac{2^{r_{-}(A)}
\cdot \sharp A^{+} [2]}{(A : A^{+} + A^{-})}
= \frac{2^{r_{-}(A)} \cdot \sharp A^{+} [2]}{(NA : 2 A^{+})} \\
&= 2^{r_{-}(A) - r_{+}(A)} \cdot (A^{+} : NA) = 2^{r(A) - 2 r_{+}(A)}
\cdot (A^{+} : NA),  \\
\text{and} \quad &h(G, A) = 2^{2 r_{+}(A) - r(A)}.
\end{align*}
In particular, if $ r_{+}(A) = 0 $ and $ A^{+}[2] = \{ O \},
$ then $ A = A^{+} + A^{-} $ and $ \sharp \text{H}^{1} (G, A) =
2^{r_{-}(A)} = 2^{r(A)}. $
\par \vskip 0.1 cm
{\bf Proof.} \ This is a consequence of the Theorem on Herbrand
quotient (see [AW]), for the detail, see the proof of
Theorem 1.5 of [Q] on a special case. \quad $ \Box $
\par  \vskip 0.2cm

{\bf Remark.} For the case that $ G $ is a group of order $ 2, $ the result
in Prop.2 is unconditional, and more explicit than the one in Prop.1 above.
For a cyclic group $ G $ of order bigger than $ 2, $ at present, by the
same way as in Prop.2 above, we only has a relatively crude result as follows,
\begin{align*} &\sharp \text{H}^{1} (G, A) = \frac{(_{N}A : (1-\sigma )(_{N}A))}
{(NA : n \cdot A^{G})}, \\
&(NA : n \cdot A^{G}) = (A : A^{G} + _{N}A) =
((1-\sigma )A : (1-\sigma )(_{N}A)),
\end{align*}
where the $ G-$module $ (A, +) $ is a finitely generated abelian
group, $ N = N_{G} \in \Z[G] $ as above,
$  _{N}A = \{a \in A: Na = 0 \}, \ NA = \{N a : a \in A \}. $
\par  \vskip 0.2cm

The proof of Proposition 1 (see Proposition 2.1 below) is given in Section 2,
and some examples are given there. In Section 3, we apply the formula in Prop.2
above to study unit groups over quadratic extensions of number fields
(see Cor.3.1 below), from which we give different proofs for some known
theorems about number fields (see the proofs of Thm.3.2 and Thm.3.5 below),
as well as some applications of Pell equations, units in CM number
fields and class number of number fields (see Prop.3.3, Prop.3.4
and Prop.3.8 below).

\par     \vskip  0.4 cm

\hspace{-0.8cm}{\bf 2. $ H^{1} $ and Herbrand quotient}

\par \vskip 0.2 cm

Let $ G $ be a group, and $ X $ be a $ G-$set. For each $ x \in X, $
the $ G-$orbit of $ x $ is $ \mathcal{O}(x) = \{ \tau x : \ \tau \in G \}, $
and the stabilizer of $ x, $ denoted by $ G_{x}, $ is the subgroup
$ G_{x} = \{ \tau \in G : \ \tau x = x\}. $ If $ X $ is finite, then
$ \sharp \mathcal{O}(x) = (G : G_{x}), $ the index of $ G_{x} $ in $ G. $
Moreover, the number $ N $ of $ G-$orbits of $ X $ is given by the following
formula: \ $ N = \frac{1}{\sharp G} \cdot \sum _{\tau \in G} F(\tau ), $
where for each $ \tau \in G, \ F(\tau ) = \sharp \{x \in X : \ \tau x = x \} $
(see [Ro, pp. 56$\thicksim$59]). In the following, we denote the set
of $ G-$orbits of $ X $ by $ X / G, $ and we write $ X / G =
\{x_{1}, \cdots, x_{r} \} $ if each $ x_{i} $ is a representative in
a $ G-$orbit, and $ x_{1}, \cdots, x_{r} $ represent all distinct
$ G-$orbits of $ X. $

Now we come to prove the Proposition 1 above, that is, the following result.
\par \vskip 0.2 cm

{\bf Proposition 2.1.} \ Let $ G $ be a finite cyclic
group, and $ A $ be a $ G-$module. If there exists a $ G-$submodule $ B $
of finite index in $ A $ which is a finitely generated free abelian group
containing a $ G-$invariant $ \Z-$basis $ X = \{x_{1}, \cdots, x_{r} \} $
(i.e., $ \tau X \subset X $ for all $ \tau \in G, $ or in other words,
$ X $ is a $ G-$set). Then
$$ h(G, A) = \prod _{x \in X / G} \sharp G_{x}, \quad \text{and} \quad
\sharp \text{H}^{1}(G, A) = \frac{(A^{G} : N_{G}A)}
{\prod _{x \in X / G} \sharp G_{x}}, $$ where $ G_{x} = \{\tau \in G :
\ \tau x = x \} $ is the stabilizer of $ x. $
\par \vskip 0.1 cm
{\bf Proof.} \ From the exact sequence of $ G-$modules $ 0 \rightarrow
 B \rightarrow A \rightarrow A / B \rightarrow 0, $ by Herbrand's theorem
(see e.g., [AW, p.109]), $ h(G, A) = h(G, B) \cdot h(G, A / B). $ By
assumption, $ A / B $ is finite, so $ h(G, A / B) = 1, $ and then
$ h(G, A) = h(G, B). $ So we only need to work out $ h(G, B). $ \\
For each $ x \in X / G, $ its $ G-$orbit is $ \mathcal{O}(x) =
\{ \tau x : \ \tau \in G / G_{x} \}. $ Let $ B(x) =
\bigoplus _{\tau \in G / G_{x}} \Z \tau x, $ then obviously, $ B(x) $
is a $ G-$submodule of $ B, $ and $ B =
\bigoplus _{x \in X / G} B(x) $ as $ G-$modules. Note that $ \Z x $ is
a trivial $ G_{x}-$module, and $ B(x) =
\bigoplus _{\tau \in G / G_{x}} \tau (\Z x) = \text{Ind}_{G_{x}}^{G}
(\Z x), $ that is, $ B(x) $ is the induced $ G-$module by the
$ G_{x}-$module $ \Z x $ (see [Br, p.67], [Neu, pp.10,11]). So by
Shapiro's lemma (see [Br, p.136], [Neu, p.11]), we have
$ \widehat{\text{H}}^{i}(G, B(x)) = \widehat{\text{H}}^{i}(G_{x}, \Z x)
(i = 0, -1). $ Hence $ h(G, B(x)) = h(G_{x}, \Z x) = \sharp G_{x}, $
the last equality follows easily from the fact that $ \Z x $ is a trivial
$ G_{x}-$module. Therefore, by Herbrand's theorem,
$$ h(G, B) = h(G, \bigoplus _{x \in X / G} B(x)) = \prod _{x \in X / G}
h(G, B(x)) = \prod _{x \in X / G} \sharp G_{x}, $$ hence
$ h(G, A) = \prod _{x \in X / G} \sharp G_{x}. $ Moreover, by
definition, $ \widehat{\text{H}}^{0}(G, A) = A^{G} / N_{G}A, $ so
$ \sharp \text{H}^{1}(G, A) = \frac{\sharp \widehat{\text{H}}^{0}(G, A)}
{h(G, A)} = \frac{(A^{G} : N_{G}A)}{\prod _{x \in X / G} \sharp G_{x}}, $
and the proof is completed. \quad $ \Box $

A question here is on what conditions does a $ G-$module satisfy the
assumption of the Proposition 2.1?
\par \vskip 0.2 cm

{\bf Example 2.2.} \ For a finite cyclic group $ G, $ let $ X $
be a finite $ G-$set. Let $ A = \oplus _{x \in X} \Z x $ be the
free abelian group with the $ \Z-$basis $ X. $ Let $ G $ acts on
$ A $ by the following way: $ \tau (\sum _{x \in X} a_{x}x)
= \sum _{x \in X} a_{x} \tau (x) \ (\forall \tau \in G, x \in X, a_{x} \in
\Z). $ Then $ A $ is a $ G-$module satisfies the assumption of the
above Prop.2.1, so $ h(G, A) = \prod _{x \in X / G} \sharp G_{x}. $
\par \vskip 0.2 cm

{\bf Example 2.3.} \ Let $ K / \Q $ be a cyclic extension of degree $ n, $
its Galois group $ G = \text{Gal}(K / \Q) = < \sigma > $ with a generator
$ \sigma . $ We denote the integral ring of $ K $ by $ O_{K}. $ Assume that
$ K $ has a normal integral basis, i.e., there is
an element $ \alpha \in O_{K} $ such that $ \{\alpha , \sigma (\alpha ),
\cdots , \sigma ^{n-1} (\alpha ) \} $ is an integral basis of $ K. $ Then
$ O_{K} $ is a $ G-$module which is a finitely generated free abelian group
containing a $ G-$invariant $ \Z-$basis $ S = \{\alpha , \sigma (\alpha ),
\cdots , \sigma ^{n-1} (\alpha ) \}. $ Obviously, $ S $ is a transitive
$ G-$set, i.e., $ \sharp (S / G) = 1, $ so $ S / G = \{\alpha \} $ and
$ G_{\alpha } = \{1 \}. $ Hence by Prop.2.1 above, the Herbrand quotient
$ h(G, O_{K}) = \sharp G_{\alpha } = 1, $ and the order of the cohomology
group $ \sharp \text{H}^{1}(G, O_{K}) = \frac{(O_{K}^{G} : N_{G}O_{K})}
{1} = (\Z : \text{Tr}_{K / \Q}O_{K}). $ Here $ N_{G} \alpha =
\sum _{\tau \in G} \tau \alpha = \text{Tr}_{K / \Q} \alpha $ is the
trace of $ \alpha $ for every $ \alpha \in O_{K}. $ \\
Such cyclic number fields with a normal integral basis include the
quadratic number fields $ K $ with odd discriminant $ d(K) $ (i.e., $ 2 \nmid
d(K)$), and the cyclotomic number fields $ \Q(\zeta _{n}), $ where $ n = p $
or $ 2 p, p $ is an odd prime number, $ \zeta _{n} $ is a primitive $ n-$th
root of unity (see [Fe, p.27], [Nar, pp.165, 166]).
\par \vskip 0.2 cm

{\bf Example 2.4.} \ Let $ A $ be an abelian variety defined over $ \Q, $
and let $ K $ be a cyclic number field with a cyclic Galois group
$ G = \text{Gal}(K / \Q). $ By Mordell-Weil theorem, the set $ A(K) $
of $ K-$rational points of $ A $ is a finitely generated abelian group.
Assume that $ A(K) $ has a free
$ \Z-$basis $ \{P_{1}, \cdots , P_{r} \} (r = \text{rank} A(K)), $ such that
$ S = \{P_{1}, \cdots , P_{r} \} $ is a $ G-$set. Then $ A(K) $ is a
$ G-$module satisfying the assumption of the above Prop.2.1. We may
as well assume that $ S / G = \{P_{1}, \cdots , P_{r_{1}} \} $ with
$ r_{1} \leq r. $ For each $ i \in \{1, \cdots , r_{1} \}, $ let
$ K_{i} = \Q (P_{i}) $ be the defined field of $ P_{i}. $ Then the
stabilizer $ G_{P_{i}} = \text{Gal}(K / K_{i}). $ So by Prop.2.1
above, the Herbrand quotient $ h(G, A(K)) =
\prod _{i = 1}^{r_{1}} \sharp G_{P_{i}} = \prod _{i = 1}^{r_{1}}
\sharp \text{Gal}(K / K_{i}) = \prod _{i = 1}^{r_{1}} [K : K_{i}], $
and $ \sharp \text{H}^{1}(G, A(K)) = \frac{(A(K)^{G} : N_{G}A(K))}
{\prod _{i = 1}^{r_{1}} [K : K_{i}]} = \frac{(A(\Q) : N_{K / \Q}A(K))}
{\prod _{i = 1}^{r_{1}} [K : K_{i}]}, $ in particular,
$ (A(\Q) : N_{K / \Q}A(K)) \geq \prod _{i = 1}^{r_{1}} [K : K_{i}]. $ \\
If all $ P_{1}, \cdots , P_{r} \in A(\Q), $ then $ \text{rank} A(K) =
\text{rank} A(\Q), $ and $ r_{1} = r, $ so $ K_{i} = \Q, G_{P_{i}} =
G $ for each $ i = 1, \cdots , r. $ Hence in this case, $ h(G, A(K))
 = [K : \Q]^{r} = n^{r}, n = [K : \Q]. $

\par     \vskip  0.4 cm

\hspace{-0.8cm}{\bf 3. $ S-$Unit groups }

\par \vskip 0.2 cm

For the quadratic extension $ K / F $ of number fields with $ K =
F(\sqrt{D}) $ for some $ D \in F^{\ast } \setminus F^{\ast^{2}}, $
let $ G= \text{Gal}(K/F) = <\sigma > = \{1, \sigma \} $ be its
Galois group with a generator $ \sigma . $ Let $ S_{F} $ be a finite
set of primes of $ F, $ always containing all infinite primes of $
F. $ Let $ S_{K} $ be the set of primes of $ K $ lying above those
in $ S_{F}. $ The group of $ S_{K}-$units of $ K $ (resp. $
S_{F}-$units of $ F $) is denoted by $ \mathbf{U}_{K, S} \ $ (resp.
$ \mathbf{U}_{F, S} $). In particular, if $ S_{F} = S_{\infty } $
consists of all infinite primes of $ F, $ then we simply write $
\mathbf{U}_{K, S_{\infty }} = \mathbf{U}_{K} $ (resp. $
\mathbf{U}_{F, S_{\infty }} = \mathbf{U}_{F} $). By Dirichlet unit
theorem (see [Neu], p.73), $ \mathbf{U}_{K, S} \ $ and $
\mathbf{U}_{F, S} $ are finitely generated abelian groups. For the
ranks, We denote $ r_{F,S}= \text{rank}(\mathbf{U}_{F,S}) $ and $
r_{K,S}= \text{rank}(\mathbf{U}_{K,S}), $ We simply write $
r_{F,S_{\infty }} = r_{F} $ and $ r_{K,S_{\infty }} = r_{K}. $ Then
$ \mathbf{U}_{F,S} \simeq \mathbf{W}_{F} \times \Z ^{r_{F,S}}, \
\mathbf{U}_{K,S} \simeq \mathbf{W}_{K} \times \Z ^{r_{K,S}}, $ where
$ \mathbf{W}_{F} $ and $ \mathbf{W}_{K} $ are the groups of roots of
unity in $ F $ and $ K, $ respectively, $ r_{F,S} = \sharp S_{F} -1,
\ r_{K,S} = \sharp S_{K} -1. $ Obviously, $ \mathbf{U}_{K,S} $ is a
$ G-$module (see [Neu], p.74), and we have $ \mathbf{U}_{K,S}^{G} =
\mathbf{U}_{F,S}. $ We also denote $ \mathbf{U}_{K,S}^{-} = \{
\alpha \in \mathbf{U}_{K,S} : \ N (\alpha ) = 1 \}, \ N
\mathbf{U}_{K,S} = \mathbf{U}_{K,S}^{1+ \sigma } = \{ N (\alpha ): \
\alpha \in \mathbf{U}_{K,S} \}, \ N^{-} \mathbf{U}_{K,S} = \{ \alpha
/ \sigma (\alpha):
 \alpha \in \mathbf{U}_{K,S} \}, $ where $ N (\alpha ) = N_{K/F}(\alpha) =
\alpha \cdot \sigma (\alpha) $ is the norm of the element $ \alpha \in K $
over $ F. $ We write $ r^{-}_{K,S}= \text{rank}(\mathbf{U}_{K,S}^{-}), \
r^{-}_{K,S_{\infty }} = r^{-}_{K}, $
we also write simply $ \mathbf{U}_{K,S_{\infty }}^{-} =
\mathbf{U}_{K}^{-}, N \mathbf{U}_{K,S_{\infty }} = N \mathbf{U}_{K},
 N^{-} \mathbf{U}_{K,S_{\infty }} = N^{-} \mathbf{U}_{K}. $
Obviously, $ \mathbf{U}_{F,S}[2] = \{\pm 1 \}, $ and $ r_{K,S} = r_{F,S}
+ r^{-}_{K,S}, $ i.e., $ r^{-}_{K,S} = \sharp S_{K} - \sharp S_{F}. $
Also $ r(N \mathbf{U}_{K}) = r_{F}. $
\par  \vskip 0.2 cm

{\bf Corollary 3.1.} \ The order of the group
$ \text{H}^{1} (G, \mathbf{U}_{K,S}) $ is
\begin{align*} \sharp \text{H}^{1} (G, \mathbf{U}_{K,S})
&= \frac{2^{\sharp S_{K} - \sharp S_{F} + 1}}{(\mathbf{U}_{K,S} :
\mathbf{U}_{F,S} \cdot \mathbf{U}_{K,S}^{-})}
= \frac{2^{\sharp S_{K} - \sharp S_{F} + 1}}{(N \mathbf{U}_{K,S} :
\mathbf{U}_{F,S}^{2})} \\
&= 2^{ r_{K,S} - 2 r_{F,S}} \cdot (\mathbf{U}_{F,S} : N \mathbf{U}_{K,S})
= 2^{ \sharp S_{K} - 2 \sharp S_{F} + 1} \cdot
(\mathbf{U}_{F,S} : N\mathbf{U}_{K,S}).
\end{align*}
In particular,  $ (\mathbf{U}_{K,S} : \mathbf{U}_{F,S} \cdot
\mathbf{U}_{K,S}^{-}) \mid 2^{\sharp S_{K} - \sharp S_{F} + 1} $ and
$ h(G, \mathbf{U}_{K,S}) = 2^{ 2 \sharp S_{F}- \sharp S_{K}  - 1}. $

{\bf Proof.} Easily follows from Dirichlet unit theorem and the
above Prop.2. \ $ \Box $
\par  \vskip 0.2 cm

For the first application, we consider the known formula in
the following Thm.3.2 for the Herbrand quotient of the relative
quadratic extension, and gives a different proof of Theorem 1.3 on
p.74 of [Neu] in this case.
\par  \vskip 0.2 cm

{\bf Theorem 3.2} (see the Theorem 1.3 on p.74 of [Neu] for $ n = [K : F] = 2 $). \
The Herbrand quotient
$$ h(G, \mathbf{U}_{K,S}) = \frac{1}{2} \prod _{v \in S_{F}} n_{v}, $$ \
where $  n_{v} $ is the order of the decomposition group
$ G_{w} \subset G $ above $ v. $
\par  \vskip 0.1 cm

{\bf Proof.} \ On the one hand,
$ S_{F} = S_{r} \sqcup S_{i} \sqcup S_{s}, $
where $ S_{r} = \{ v \in S_{F} : v \ \text{ramifies in} \ K \}, \
S_{i} = \{ v \in S_{F} : v \ \text{is inertia in} \ K \}, S_{s} =
\{ v \in S_{F} : v \ \text{splits in} \ K \}. $ So $ \sharp S_{F} =
\sharp S_{r} + \sharp S_{i} + \sharp S_{s} $ and $ \sharp S_{K} =
\sharp S_{r} + \sharp S_{i} + 2 \cdot \sharp S_{s}. $ Also, $ n_{v} =
1 $ for $ v \in S_{s}, $ and $ n_{v} = 2 $ for $ v \in S_{r} \cup S_{i}, $
so $ \prod _{v \in S_{F}} n_{v} = 2^{\sharp S_{r} + \sharp S_{i}} =
2^{2 \cdot \sharp S_{F} - \sharp S_{K}}. $ On the other hand, by Cor.3.1
above, $ h(G, \mathbf{U}_{K,S}) = 2^{ 2 \sharp S_{F}- \sharp S_{K}  - 1}, $
so the equality holds.
\quad $ \Box $
\par  \vskip 0.2 cm

For the second application, we consider the case of real quadratic
number fields.

{\bf Proposition 3.3.} \ Let $ K = \Q (\sqrt{D}) $ be a real
quadratic number field with a square-free integer $ D (> 0), $ and $
G = \text{Gal}(K/\Q). $ Let $ \epsilon > 1 $ be a fundamental unit
of $ K. $ Then the following statements are equivalent: \\
(1) \ $ N(\varepsilon) = -1; $ \\
(2) \ $ \sharp \text{H}^{1} (G, \mathbf{U}_{K}) = 2; $ \\
(3) \ The Pell equation $ x^{2} - y^{2} D = -1 $ in the case $ D
\equiv 2, 3 \ (\text{mod} \ 4) $ \ (respectively, $ (2x - y)^{2} -
y^{2} D = -4 $ in the case $ D \equiv 1 \ (\text{mod} \ 4 $)) has
integral solutions.

\par  \vskip 0.1cm
{\bf Proof.} \ By Dirichlet unit theorem, $ \mathbf{U}_{K} = <-1>
\times <\epsilon >, $ this $\epsilon $ is the unique smallest unit
greater than $1$ (see [Wei], pp.238,239 and [JW], pp.82,83).
Obviously, $ N(\varepsilon) = \pm 1, $ so $ N\mathbf{U}_{K} \subset
\{ \pm 1\}. $ It follows from the formula in the above Cor.3.1 that
$ \sharp \text{H}^{1} (G, \mathbf{U}_{K}) = 2 (\mathbf{W}_{\Q} :
N\mathbf{U}_{K}) $ with $ \mathbf{W}_{\Q} = \{ \pm 1 \}. $ So $
\sharp \text{H}^{1} (G, \mathbf{U}_{K}) = 2 $ if and only if $
N\mathbf{U}_{K} = \{ \pm 1 \}. \\
$ (1) $ \Leftrightarrow $ (2). If there is an element $ \alpha \in
\mathbf{U}_{K} $ such that $ N \alpha = -1, $ then since $ \alpha =
\pm \epsilon ^{m} $ for some $ m \in \Z, $ we have $ -1 = N(\epsilon
)^{m}, $ so $ m $ is odd and $ N(\epsilon ) = -1. $ It follows from
this that $ N\mathbf{U}_{K} = \{ \pm 1\} $ if and only if $
N(\epsilon ) = -1. $ Hence $ \sharp \text{H}^{1} (G, \mathbf{U}_{K})
= 2 $ if and only if $ N(\epsilon ) = -1. $ \\
(2) $ \Leftrightarrow $ (3). If $ D \equiv 1 \ (\text{mod} \ 4), $
then the integral ring $ \mathcal{O}_{K} = \Z[\frac{-1 + \sqrt{D}
}{2}], $ each element of $ \mathcal{O}_{K} $ has the form $ x + y
\cdot \frac{-1 + \sqrt{D} }{2}, x, y \in \Z. $ So $ N\mathbf{U}_{K}
= \{ \pm 1 \} $ if and only if there are $ x, y \in \Z $ such that $
N( x + y \cdot \frac{-1 + \sqrt{D} }{2}) = -1, $ i.e., the equation
$ (2x - y)^{2} - y^{2} D = -4 $ has integral solutions. Similar for
the case $ D \equiv 2, 3 \ (\text{mod} \ 4), $ and the conclusion
follows from the above discussion.  \quad $ \Box $ \\
Note that the set of $ D > 0 $ for which the norm of the fundamental
unit $ \epsilon $ is $ -1 $ has not been determined (see [IR],
P.192). Many interesting facts about Pell equation and its
applications can be found in [JW].
\par  \vskip 0.2cm

For the third application, we consider the case of CM number fields.

{\bf Proposition 3.4.} \ Let $ K $ be a CM number field (see [Wa],
p.38), $ F = K^{+} $ the real subfield, $ [K : F] = 2, $ and the
Galois group $ G = \text{Gal}(K/F) = <\sigma > $ is generated by the
complex conjugation $ \sigma $ on $ K. $ If $ \sharp S_{F} = \sharp
S_{K}, $ for example, when every finite prime in $ S_{F} $ is
ramified or inertia in $ K. $ Then
$$ \sharp \text{H}^{1} (G, \mathbf{U}_{K,S}) = \frac{2}
{(\mathbf{U}_{K,S} : \mathbf{U}_{F,S} \cdot \mathbf{W}_{K})} =
\frac{2} {(N \mathbf{U}_{K,S} : \mathbf{U}_{F,S}^{2})} =
2^{-r_{K,S}} \cdot (\mathbf{U}_{F,S} : N \mathbf{U}_{K,S}), $$ in
particular, \ $ (\mathbf{U}_{K,S} : \mathbf{U}_{F,S} \cdot
\mathbf{W}_{K}) = (N \mathbf{U}_{K,S} : \mathbf{U}_{F,S}^{2}) \mid
2, \ \sharp \text{H}^{1} (G, \mathbf{U}_{K,S}) \mid 2, $ and $
(\mathbf{U}_{F,S} : N \mathbf{U}_{K,S}) = 2^{r_{K,S}} $ or $
2^{r_{K,S} + 1}. $
\par  \vskip 0.1cm
{\bf Proof.} \ Since $ \text{rank}(\mathbf{U}_{K,S}^{-}) = r_{K,S} -
r_{F,S} = \sharp S_{K} - \sharp S_{F} = 0 $ and $ F = K^{+} $ is
totally real, it is easy to see that $ \mathbf{U}_{K,S}^{-} =
\mathbf{W}_{K}. $ Note that $ \sharp \text{H}^{1} (G,
\mathbf{U}_{K,S}) $ is a positive integer number, the conclusion
follows directly from the formula of the above Cor.3.1. \quad $ \Box
$
\par  \vskip 0.2cm

From the above Prop.3.4, one can directly deduces the known index estimation of
the unit groups of $ F, K, $ which gives a different proof of
Theorem 4.12 on p.40 of [Wa].
\par  \vskip 0.2cm

{\bf Theorem 3.5} (see the Theorem 4.12 on p.40 of [Wa]). \ Let $ K
$ and $ F $ be as in Prop 3.4 above. Then $ (\mathbf{U}_{K} :
\mathbf{U}_{F} \cdot \mathbf{W}_{K}) \mid 2. $
\par  \vskip 0.1 cm

{\bf Proof.} \ Take $ S_{F} = S_{\infty }, $ then $
\mathbf{U}_{F,S} = \mathbf{U}_{F}, \mathbf{U}_{K,S} =
\mathbf{U}_{K}, $ and the condition of the above Prop.3.4 holds, so
the conclusion follows. \quad $ \Box $
\par  \vskip 0.2 cm

{\bf Remark.} For the cyclotomic field case, i.e., $ K = \Q(\zeta
_{m}) \ ( 2 < m \in \Z) $ and $ F = K^{+} = \Q (\zeta _{m} + \zeta
_{m}^{-1}), $ where $ \zeta _{m} $ is a primitive $ m-$th root of
unity in the complex number field $ \C. $ If $ m = p^{t} $ for some
prime number $ p $ and integer $ t \geq 1, $ then $ \mathbf{U}_{K} =
\mathbf{U}_{F} \cdot \mathbf{W}_{K} $ (see [Wei], p.268 and [Wa],
p.40). So by the above Prop.3.4 we have $ \sharp \text{H}^{1} (G,
\mathbf{U}_{K}) = 2, \
 N \mathbf{U}_{K} =  \mathbf{U}_{F}^{2},  $ and
 $ (\mathbf{U}_{F} : N \mathbf{U}_{K}) = 2^{r_{K} + 1}. $
\par \vskip 0.2 cm

Now for the CM field $ K $ and its real subfield $ F = K^{+} $ as
above, let $ \mathbf{U}_{F}^{+} $ be the set consisting of all
totally positive units in $ F. $ It is easy to see that $
\mathbf{U}_{F}^{2} \subset N \mathbf{U}_{K} \subset
\mathbf{U}_{F}^{+} \subset \mathbf{U}_{F}. $ Let $ \phi : \
\mathbf{U}_{K} \rightarrow \mathbf{W}_{K}, \ \alpha \mapsto \alpha /
\sigma (\alpha ). $ Then from the proof of Theorem 4.12 in [Wa,
p.40], we have
\par \vskip 0.2 cm

{\bf Corollary 3.6.} \ (1) \ If $ \phi (\mathbf{U}_{K}) =
\mathbf{W}_{K}, $ then $ (\mathbf{U}_{F} : N \mathbf{U}_{K}) =
2^{r_{K}}, \ (N \mathbf{U}_{K} : \mathbf{U}_{F}^{2})=
(\mathbf{U}_{K} : \mathbf{U}_{F} \cdot \mathbf{W}_{K}) = 2, $ and $
\sharp \text{H}^{1} (G, \mathbf{U}_{K}) = 1.  $ In particular,
$ \mathbf{U}_{F}^{+} \neq \mathbf{U}_{F}^{2}; $ \\
(2) \ If $ \phi (\mathbf{U}_{K}) = \mathbf{W}_{K}^{2}, $ then $
(\mathbf{U}_{F} : N \mathbf{U}_{K}) = 2^{r_{K}+1}, \ N
\mathbf{U}_{K} = \mathbf{U}_{F}^{2}, \ \mathbf{U}_{K} =
\mathbf{U}_{F} \cdot \mathbf{W}_{K}, $ and $ \sharp \text{H}^{1} (G,
\mathbf{U}_{K}) = 2.  $
\par \vskip 0.2 cm

For example, let $ K = \Q(\zeta _{m}) $
and $ F = K^{+} = \Q (\zeta _{m} + \zeta _{m}^{-1}) $ be as above. \\
(a) \ if $ m $ is not a prime power, then the conclusion (1) of Cor.3.6 holds. \\
(b) \ if $ m $ is a prime power, then the conclusion (2) of Cor.3.6 holds.
\par  \vskip 0.2cm

{\bf Example 3.7.} \ In this example, we let $ K, F $ be two real
Galois extensions over $ \Q, \ [K : F] = 2, $ and $ [F : \Q] = n. $
The Galois group $ G= \text{Gal}(K/F) = \{1, \sigma \} $ with a
generator $ \sigma . $ Let $ \mathbf{U}_{K}, \mathbf{U}_{K}^{1},
\mathbf{U}_{K}^{+} $ and $ \mathbf{U}_{K}^{2} $ be the unit group,
the subgroup of the units of norm $ 1, $ the subgroup of totally
positive units and the subgroup of the unit squares, respectively,
of $ K. $ Similar for the meaning of $ \mathbf{U}_{F},
\mathbf{U}_{F}^{1}, \mathbf{U}_{F}^{+} $ and $ \mathbf{U}_{F}^{2}. $
Let $ \mathbf{U}_{K}^{+1} = \mathbf{U}_{K}^{+} \cap
\mathbf{U}_{K}^{1}, \ \mathbf{U}_{K}^{1, 2} = \mathbf{U}_{K}^{1}
\cap \mathbf{U}_{K}^{2}, \ r_{K}^{+1} =
\text{rank}(\mathbf{U}_{K}^{+1}), \ r_{K}^{1, 2} =
\text{rank}(\mathbf{U}_{K}^{1, 2}), \ r_{K}^{1} =
\text{rank}(\mathbf{U}_{K}^{1}), \ r_{F}^{1} =
\text{rank}(\mathbf{U}_{F}^{1}), \ r_{K}^{+} =
\text{rank}(\mathbf{U}_{K}^{+}), \ r_{F}^{+} =
\text{rank}(\mathbf{U}_{F}^{+}). $ Then we have
\begin{align*} & \sharp \text{H}^{1} (G, \mathbf{U}_{K}^{+}) = \frac{2^{n}}
{(\mathbf{U}_{K}^{+} : \mathbf{U}_{F}^{+} \cdot
\mathbf{U}_{K}^{+1})} = \frac{2^{n}} {(N \mathbf{U}_{K}^{+} :
\mathbf{U}_{F}^{+2})} =
2 \cdot (\mathbf{U}_{F}^{+} : N \mathbf{U}_{K}^{+}), \\
& \sharp \text{H}^{1} (G, \mathbf{U}_{K}^{2}) = \frac{2^{n}}
{(\mathbf{U}_{K}^{2} : (\mathbf{U}_{K}^{2})^{G} \cdot \mathbf{U}_{K}^{1, 2})} =
\frac{2^{n}}
{((N \mathbf{U}_{K})^{2} : ((\mathbf{U}_{K}^{2})^{G})^{2})} =
2 \cdot ((\mathbf{U}_{K}^{2})^{G} : (N \mathbf{U}_{K})^{2}), \\
& \sharp \text{H}^{1} (G, \mathbf{U}_{K}^{1}) = \frac{2^{n+1}}
{(\mathbf{U}_{K}^{1} : \mathbf{U}_{F} \cdot \mathbf{U}_{K}^{+1})} =
\frac{2^{n+1}}
{(N \mathbf{U}_{K}^{1} : \mathbf{U}_{F}^{2})} =
2 \cdot (\mathbf{U}_{F} : N \mathbf{U}_{K}^{1}).
\end{align*}
\par  \vskip 0.1cm
{\bf Proof.} \ By assumption, $ r_{F} = n-1, r_{K} = 2n -1. $ It is
easy to see that $ \mathbf{U}_{K}^{+}, \mathbf{U}_{K}^{1} $ and $
\mathbf{U}_{K}^{2} $ are $ G-$modules, and $
(\mathbf{U}_{K}^{+})^{G} = \mathbf{U}_{F}^{+}, \
(\mathbf{U}_{K}^{1})^{G} = \mathbf{U}_{F}. $ For $ \alpha \in
\mathbf{U}_{K}, $ the norm $ N_{K/\Q}(\alpha ) = \pm 1 $ (see [Nar],
p.96), so $ \mathbf{U}_{K}^{2} \subset \mathbf{U}_{K}^{1}. $ The
following facts are also easy: $ r_{K}^{+} = r_{F}^{+} + r_{K}^{+1},
\ r_{K}^{1, 2} = r_{K}^{+1} = r_{K} - r_{F} = n, \ r_{K}^{1} = r_{K}
= 2n -1, \ r_{F}^{1} = r_{F} = n-1, \ \mathbf{U}_{F}^{+} \cap
\mathbf{U}_{K}^{+1} = \mathbf{U}_{K}^{+}[2] = \{1 \}, \
\mathbf{U}_{F} \cap \mathbf{U}_{K}^{+1} =
\mathbf{U}_{F}[2]=\mathbf{U}_{K}^{+1}[2] = \{1, -1 \}. $ Then the
conclusions follow directly from the above Prop.2.   \quad $ \Box
$
\par \vskip 0.2 cm

Now again let $ F, K, G, S_{F}, S_{K}, \mathbf{U}_{F, S},
\mathbf{U}_{K, S} $ be as in the above Cor.3.1. \\
Recall that $ S_{F} $ is called large for $ K/F $ if \\
 (i) \ $ S_{F} $ contains all ramified primes of $ K/F; $ \\
 (ii) \ the $ S_{K}-$class group of $ K $ is trivial. \\
Moreover, if $ S_{F} $ is large for $ K/F, $ and satisfies the following condition \\
(iii) \ $ G = \cup _{w \in S_{K}} G_{w}, $ where $ G_{w} $ is the decomposition
group at the place $ w, $ \\
then $ S_{F} $ is called larger for $ K/F. $ (See [We,
pp.1$\sim$2]). We denote the $ S_{F}-$class group of $ F $ by $
Cl_{F,S}, $ which is defined to be the quotient of the ideal class
group of $ F $ by the subgroup generated by the classes represented
by the finite primes in $ S_{F} $ (see [Sa], p.127).
\par  \vskip 0.2cm

{\bf Proposition 3.8.} \ If $ S_{F} $ is large for $ K/F, $ then the
$ S_{F}-$class number of $ F $
$$ h_{F,S} = \sharp Cl_{F,S} = \sharp \text{H}^{1} (G, \mathbf{U}_{K,S}) =
2^{ \sharp S_{K} - 2 \sharp S_{F} + 1} \cdot (\mathbf{U}_{F,S} :
N\mathbf{U}_{K,S}), $$ so $  h_{F,S} $ is a $ 2-$power. In
particular, if $ K/F $ is ramified at some finite place, and $ S_{F}
$ is large, then
$$ (\mathbf{U}_{F,S} : N\mathbf{U}_{K,S}) =
2^{2 \sharp S_{F} - \sharp S_{K} - 1}. $$
\par  \vskip 0.1cm
{\bf Proof.} \ By Lemma 1 of [We, p.9], we have $ \text{H}^{1} (G,
\mathbf{U}_{K,S}) \simeq Cl_{F,S}, $ so by Cor.3.1 above, the first
formula follows. For the second formula, by our assumption, it is
obvious that $ S_{F} $ is also larger for $ K/F. $ By Lemma 1 of
[We, p.9], the $ S_{F}-$class group of $ F $ is trivial, so $ \sharp
\text{H}^{1} (G, \mathbf{U}_{K,S}) = 1, $ which implies that $
(\mathbf{U}_{F,S} : N\mathbf{U}_{K,S}) = 2^{2 \sharp S_{F} - \sharp
S_{K} - 1}. $ \quad $ \Box $

\par  \vskip 0.3 cm

\hspace{-0.8cm} {\bf References }
\begin{description}

\item[[AW]] M. Atiyah, C.T.C. Wall, Cohomology of groups, in:
Algebraic Number Theory (J.W.S. Cassels and A. Frohlich, Eds.),
pp.94-115, London: Academic Press, 1967.

\item[[Br]] K.S. Brown, Cohomology of Groups, New York: Springer-Verlag,
1982.

\item[[Fe]] K.Q. Feng, Algebraic Number Theory (in Chinese),
Beijing: Science Press, 2000.

\item[[IR]] K. Ireland, M. Rosen, A Classical Introduction to Modern
Number Theory, 2nd Edition, New York: Springer-Verlag, 1990.

\item[[JW]] M.J.Jacobson,Jr, H.C.Williams, Solving the Pell
Equation, CMS Books in Math., Springer, 2009.

\item[[Nar]] W. Narkiewicz, Elementary and Analytic Theory of
Algebraic Numbers, 3rd Edition, New York: Springer-Verlag, 2004.

\item[[Neu]] J. Neukirch, Class Field Theory, New York: Springer-Verlag, 1986.

\item[[Q]] Derong Qiu, On quadratic twists of elliptic curves and some
applications of a refined version of Yu's formula, Communications in Algebra,
42(12), 5050-5064, 2014.

\item[[Ro]] J.J. Rotman, An Introduction to the Theory of Groups,
Fourth Edition, New York: Springer-Verlag, 1995.

\item[[Sa]] J. W. Sands,  Popescu's conjecture in multi-quadratic
extensions, Contemporary Math., vol.358 (2004), 127-141.

\item[[Se]] J. -P. Serre, Local fields, New York: Springer-Verlag,
1979.

\item[[Wa]] L.C.Washington, Introduction to Cyclotomic Fields, 2nd Edition,
New York: Springer-Verlag, 1997.

\item[[We]] A. Weiss, Multiplicative Galois Module Structure, AMS, Providence,
Rhode Island, 1996.

\item[[Wei]] E. Weiss, Algebraic Number Theory, New York:
McGraw-Hill Book C ompany, Inc, 1963.

\end{description}

\end{document}